\documentclass{article}
\usepackage{amssymb}
\usepackage{amsmath}
\usepackage{amsfonts,amsthm}
\newtheorem{theorem}{Theorem}
\newtheorem{lemma}[theorem]{Lemma}

\theoremstyle{remark}
\newtheorem{remark}[theorem]{Remark}

\newtheorem{problem}{Problem}

\theoremstyle{definition}
\newtheorem{definition}[theorem]{Definition}

\newcommand{\Z}{\mathbb{Z}}
\newcommand{\Q}{\mathbb{Q}}
\newcommand{\C}{\mathbb{C}}

\begin{document}
\title{Generation of polycyclic groups}
\author{Martin Kassabov\footnote{The author was partially supported by
AMS Centennial Fellowship and NSF grant DMS~0600244.}\ \
\& Nikolay Nikolov}
\maketitle

\begin{center} To Dan Segal on occasion of his 60-th birthday
\end{center}
\begin{abstract} We give a new and self-contained proof of a theorem of Linnell
and Warhurst that $d(G) - d(\widehat
G) \leq 1$ for virtually polycyclic groups $G$. We also give a simple
sufficient condition for equality $d(G)=d(\widehat G)$ when $G$ is virtually abelian.
\end{abstract}
\section*{Introduction}

Let $G$ be a finitely generated residually finite group. By $d(G)$ we denote
the minimal size of a generating set for $G$, and by $d(\widehat G)$
the minimal size of a generating set for the profinite completion $\widehat G$ of $G$. In other words
\[
d(\widehat G)= \max \left\{ d(G/N) \ | \ N \vartriangleleft G, \ G/N
< \infty \right\}.
\]

Polycyclic groups are one of the best understood class of groups. For example
most of the decision problems are decidable in this class, see \cite{Danpaper}.

It seems surprising therefore that it is still an open problem whether there exists an algorithm which finds $d(G)$ for any polycyclic group $G$ (given by say a set of generators and relations).
This is unknown even in
the case when $G$ is virtually abelian.

It is obvious that $d(G) \geq d(\widehat G)$ and when there is equality
both the value of $d(G)$ and a minimal generating set for $G$ can indeed be found algorithmically. (Say by enumerating both
the finite images and all possibilities for generating sets for $G$).

In general $d(G)-d(\widehat
G)$ can be arbitrarily large even for metabelian  groups $G$, see \cite{Nos}.
In fact Wise~\cite{Wise} has proved that there
exist groups $G$ with arbitrarily large $d(G)$ while $d(\widehat G)=3$.

Fortunately for polycyclic groups the situation is not that bad.
In~\cite{LW} Linnell and Warhurst proved the following theorem using methods from commutative
algebra and lattices over orders.

\begin{theorem}
\label{main}
Let $G$ be a virtually polycyclic group. Then
$d(G) \leq d(\widehat G) +1$.
\end{theorem}

Note that the inequality is sharp even for virtually abelian groups: many examples with $d(G)=d(\widehat G)+1$ are constructed in \cite{Nos2}.

In this note we give an alternative proof of Theorem \ref{main}. While not claiming anything new
we believe that our argument is much simpler
that the original one in \cite{LW}. Moreover our result gives some sufficient
condition when $d(G)=d(\widehat  G)$ which can be verified quite easily in
the case when $G$ is virtually abelian.

\begin{theorem}
\label{vabelian}
Let $G$ be a group with normal finitely generated abelian subgroup $U$ such
that $G/U$ is finite. Let $d_p(G)=d(G/U^p)$ for any prime $p$.
Then
\begin{equation}
\label{abeq}
\alpha = d(\widehat G) \leq d(G) \leq
k:= \max \left\{ \alpha, \beta +1\right\}
\leq d(\widehat G) +1,
\end{equation}

where $\alpha=\max_p d_p(G)$ and $\beta  =\min_p d_p(G)$.

Moreover for any integer $N \in \mathbb N$ there exists a generating set $S$ for $G$  of size
$k$, such that the first $d(\widehat G)$ elements generate a
subgroup of index co-prime to $N$.
The same result holds for finitely generated
virtually nilpotent groups.
\end{theorem}
In particular if there are two primes $p$ and $q$ such that $d_p(G) \not
= d_q(G)$ then $d(G)= d(\widehat G)$.

Note that Theorem \ref{vabelian} easily implies a weaker version of
Theorem \ref{main}, namely that $d(G) \leq d(\widehat G) +2$, however obtaining
the right bound $d(\widehat G)+1$ is harder. For that we need a general and somewhat technical result (Theorem \ref{lifts} on lifting generators) proved
in Section \ref{lif}. 
The proofs of Theorems \ref{vabelian} and \ref{main} are then immediate and are given in Section \ref{pf}. 
\subsubsection*{Notation}
For elements
$a,b \in G$ in a group $G$ the commutator $[a,b]$ of $a$ and $b$ is $aba^{-1}b^{-1}$.

\section{Lifting generators} \label{lif}

In this section we shall prove a general result which under certain
condition produces a generating set of a group $G$ starting from a generating set of some quotient $G/V$ of $G$. The reason for stating it in such generality is because in the next section we shall apply it in two settings: when $G$ a virtually abelian group and then when $G$ is virtually metabelian.

First recall the following result by Gasch\"{u}tz (\cite{Ga}).
\begin{theorem} \label{gas}
Let $G$ be a finite group with a normal subgroup $N$.
Let $d \geq d(G)$ and let
$a_1, \ldots a_d$ be any $d$ elements which generate $G$ mod $N$, i.e.,
$G=N \langle a_1, \ldots ,a_d\rangle$.
Then we can find elements $g_i \in a_iN$ ($i=1,2, \ldots, d$)
such that $G=\langle g_1, \ldots ,g_d\rangle$.
\end{theorem}
\medskip

\begin{definition} Let $G$ be a group and $p$ be a prime. A normal subgroup $L$ of finite index
in $G$ is \emph{p-good} if for any subgroup $H \leq G$ with $HL=G$ we have
that $[G:H]$ is finite and coprime to $p$.
\end{definition}
It is not difficult to see that $p$-good subgroups exists for any virtually
polycyclic group $G$, see Lemma \ref{pgood} below.
\medskip

Now let $G$ be a finitely generated group with an abelian normal subgroup $V$ of finite rank.
Suppose that for every prime number $p$ we have chosen a $p$-good subgroup $G_p$ of $G$ such that $G_p  \geq V^p$.

\begin{definition}
We say that $h_1,\dots, h_k$ generate $G$ mod $p$ if
$\langle h_1,\dots, h_k\rangle G_p=G$. Let $d_p(G)=d(G/G_p)$ denote the minimal
size of a set of generators for $G$ mod $p$.
\end{definition}

\begin{definition}
\label{fox}
Let $w=w(x_1, \ldots , x_n)$ be a group word (element in the free group $F$).
The Fox derivatives $\frac{\partial w}{\partial x_i}$
are elements in the group ring $\Z[F]$,
which are defined by
$\frac{\partial x_j}{\partial x_i}=\delta_{ij}$
and
$$
\frac{\partial uv}{\partial x_i}=
\frac{\partial u}{\partial x_i}+
u \frac{\partial v}{\partial x_i}.
$$
Let $G$ be a group and $V$ be a $G$-module.
For any $n$-tuple $\underline{g}=(g_1, \ldots ,g_n)\in G^n$
the Fox derivative $\frac{\partial w}{\partial x_i}$ naturally defines a map
$\frac{\partial w}{\partial x_i} (\underline g): V \rightarrow V$.
\end{definition}

An equivalent way to define this map is the following: Let $\Gamma$ be an
extension of $G$ by the abelian group $V$ then
$$
\frac{\partial w}{\partial x_i}(\underline{g})(a)=
w(\gamma_1, \ldots , \gamma_{i-1},\
a \gamma_i,\gamma_{i+1},
\ldots, \gamma_n)  \cdot w(\gamma_1, \ldots, \gamma_n)^{-1}
$$
for any lifts $\gamma_i\in \Gamma$ of $g_i\in G$.

\begin{theorem}
\label{lifts}
Let $\underline{\gamma} =(\gamma_1, \dots \gamma_k)$ be a set of elements in
$G$ which generate $G/V$. Suppose that $w(x_1,\dots, x_k)$ is a word such
that $w(\underline{\gamma}) \in V$. 

Assume that
\begin{enumerate}
\item The image of the map $\pi:V \to V$ defined by
$
\pi(v)= \frac{\partial w}{\partial x_k}(\underline{\gamma}) \circ v
$
has finite index $M$ in $V$.

\item For any choice of elements $g_1 \in \gamma_1 V,\ldots, g_{k-1} \in\gamma_{k-1}V$
the group \\ $\langle g_1, \ldots ,g_{k-1} \rangle$ generated by them has
finite index in $G$.

\item We have that $d_p(G) \leq k$
for any prime $p$.
\end{enumerate}

Then there exist lifts $g_1,\dots,g_k$ of $\gamma_1, \ldots , \gamma_k$ (i.e.
such that  $g _i \in \gamma_i V$) which generate $G$.

Moreover there is an algorithm for finding $g_i$ from the $\gamma_i$ (provided
all the objects from conditions 1,2,3 above are computable).
\end{theorem}

\begin{proof}
We say that the element $g_i$ is a lift of $\gamma_i$ whenever $g_i
\in \gamma_i V$. Note that the Fox derivative
$\frac{\partial w}{\partial x_i} (\underline{g})=
\frac{\partial w}{\partial x_i} (\underline{\gamma})$
does not depend on the choice of lifts
$\underline{g}=(g_1,\ldots, g_k)$ of $\underline{\gamma}=(\gamma_1, \ldots ,\gamma_k)$.
\medskip

For each $i=0,1, \ldots, k$ let $Q(i)$ be the following statement.
\bigskip

$\mathbf{Q(i):}$ {\em There exist lifts $S_i=\{g_1, \ldots, g_i\}$ of $\gamma_1, \ldots , \gamma_i$ and a finite set of prime numbers $P_i$
with the following property.
\begin{itemize}
\item
For each prime $p\in P_i$ there exist lifts $g_j^{(p)} \in \gamma_j V$, $j=i+1,\dots, k$
such that
$S_i \cup \{g_{j}^{(p)}\}_{j>i}$ are $k$ generators for $G$ mod $p$.

\item
For each prime $p\not\in P_i$,
there exist lifts $g_j^{(p)}$, $j=i+1,\dots,k-1$
such that for any lift $g_k^{(p)} \in \gamma_k V$ we have that
$S_i \cup \{g_{j}^{(p)}\}_{j>i}$ are $k$ generators for $G$ mod $p$.
\end{itemize}}
\bigskip

The proof of $Q(i)$ is by induction on $i$.

The base case $i=0$ is proved as follows: Take $S_0= \emptyset$ and 
choose any lifts $\gamma_1,\dots,\gamma_{k-1}$ of the $\gamma_j$. They generate a subgroup of
finite index $L$ in $G$,  therefore for any $p \not | L$ and any lift
$g_k \in \gamma_k V$ the elements $\{g_1, \ldots ,g_k\}$ generate $G$ mod $p$.

Define $P_0$ to be the set consisting
of all primes which divide $L$ or $M$.
We only have to show that for all $p \in P_0$ there exist lifts $g_1^{(p)},
\ldots ,g_k^{(p)}$, which generate $G$ mod $p$. Consider the $k$ images $\bar \gamma_j$ of $\gamma_1, \ldots , \gamma_k$
in $G/VG_p$. They generate $G/VG_p$ (since $\gamma_j$ generate $G/V$) and
also we know that $G/G_p$ is $k$-generated. By Gaschutz Theorem we can find
elements $g_j^{(p)} \in \gamma_j VG_p$ which generate $G/G_p$. These can
be further adjusted by elements from $G_p$ so that $g_j^{(p)} \gamma_j^{-1}
\in V$. This proves the base case $i=0$ of the induction.
\medskip

Suppose that we have already found $S_i$ and $P_i$.
By the Chinese Remainder Theorem there exists a lift $g_{i+1}$ such that
$g_{i+1} = g_{i+1}^{(p)} (\mathrm{ mod \ }G_p) $ for all $p \in P_i$.

Choose any lifts $g_{i+2},\dots g_{k-1}$. By one of the assumptions
the group $H$ generated by
$\{g_i\}_{i=1}^{k-1}$ is of finite index $N_i$
in $G$. Denote $P_{i+1}=\{p \mid \,\, p \textrm{ divides } N_i\} \cup P_i $.

We want to show that $S_{i+1} = S_i \cup \{g_{i+1}\}$ and the set $P_{i+1}$
satisfy the induction hypothesis.

It is very easy to check that the second part of the induction
hypothesis is satisfied for this definition of the set $P_{i+1}$ (just choose
$g_j^{(p)}=g_j$ for $j=i+2, \ldots, k-1$). It remains to show that for all primes $p \in P_{i+1}$
the first condition is satisfied. This is clearly the case if $p \in P_i$.
Let $p\not\in P_i$ then by the induction assumption there exist lifts $\{g_j^{(p)}\}_{j=i}^{k-1}$
which together  with $S_i$ and any lift $g_k^{(p)} \in \gamma_k V$ generate $G$ mod $p$.
We will show that we can chose a lift $g_k$ of $\gamma_k$such that the group $L$ generated by
$S_{i+1}$, $\{g_j^{(p)}\}_{j=i+1}^{k-1}$ and $g_k$ contains an element
$u \equiv g_{i+1}^{-1} g_{i+1}^{(p)} (\mathrm{mod \ } V^p)$,
which implies that these elements generate $G$ mod $p$.

The key observation here is that for $x \in V$ the element
$w\big(g_1^{(p)},\dots, g_{k-1}^{(p)},x g_k^{(p)}\big)$ is equal to
$$
\left(\frac{\partial w}{\partial x_k}(\underline{\gamma}) \cdot x \right)
w\big(g_1^{(p)},\dots, g_{k-1}^{(p)},g_k^{(p)}\big)=\pi' (x) \in V
$$
which by one of the assumptions takes any value in $V/V^p$ as $x$ ranges over
$V/V^p$
(since $p$ does not divide the index $M$ of the image of
$\frac{\partial w}{\partial x_k}(\underline{\gamma})$ in $V$).
So we can indeed find $x \in V$ such that the element
$\pi'(x) \in V$ satisfies $\pi'(x)\equiv
g_{i+1}^{-1} g_{i+1}^{(p)}$ mod $V^p$.

This shows that there exist lifts $g_k$ which generate $G$ mod $p$
which completes the induction step.

\medskip

The statement $Q(k)$ gives a set $S_k$ which generates $G$ mod $p$ for any prime $p$ and therefore $\langle S_k \rangle = G$. 

It is clear that this argument in fact produces an algorithm for finding the set
$S_k$ in a very efficient way,
of course provided the various subgroup indices, words and maps $G \rightarrow
G/G_p$ involved in the induction are computable. 
Theorem \ref{lifts} is
proved.
\end{proof}

\begin{remark}
A slight modification of the proof gives that for any finite set
$P$ of primes such that $d_p(G) < k$, we can find lifts $g_1,\dots,g_{k-1}$
which together with any lift of $\gamma_k$
generate a subgroup of index not divisible by any prime in $P$.
\end{remark}

\begin{remark}
If $\gamma_k=e$ then we can take the word $w=x_k$. Its Fox derivative is
$\frac{\partial w}{\partial x_k} = e$  and defines the identity map from $V$ to $V$,
 which is clearly surjective.
\end{remark}

\begin{remark}
If we replace the assumption that $V$ is abelian with
$V$ nilpotent, then all results remain valid,
since a set generates a nilpotent group if and only if it generates the
abelianization of the group.
\end{remark}

\section{Applications of Theorem \ref{lifts}} \label{pf}

\subsection{Proof of Theorem \ref{vabelian}}
\begin{proof}
Clearly $d(\widehat G)= \alpha$.
Let $k:=\max \{\beta+1, \alpha\}$.

Take $V:=U^q$ where $q$ is a prime such that $\beta =d_q(G)$.
Then $V$ is $q$-good, i.e., any collection of elements which generates $G/V$ generates
a subgroup of finite index (coprime to $q$) in $G$.

Take elements $\gamma_1,
\ldots ,\gamma_{\beta}$ which generate $G/V$. Set $\gamma_i=1$ for $j=\beta+1,\ldots
, k$.
It is easy to see that the group $G$, subgroups $V$, $G_p=V^p$ (for any prime
$p$), the elements $\gamma_i$ above, and the word $w=x_k$ satisfy the conditions of Theorem \ref{lifts}. We conclude that $G$ can be generated by some lifts of $\gamma_1,
\ldots , \gamma_k$ and so $d(G) \leq k$ as claimed.
\medskip

For the second part of the theorem start with $V=U^N$ instead and with
any generating
set $\gamma_1, \ldots , \gamma_s$ for $G/V$ and again take $\gamma_j=1, \ s<j\leq
k$. The rest of the argument is similar.
\end{proof}

\subsection{Proof of Theorem \ref{main}}
We begin with the following straightforward

\begin{lemma} \label{pgood} If $G$ is a virtually polycyclic group and $p$
is a prime then $G$ has a $p$-good subgroup $L$.
\end{lemma}
\begin{proof} We use induction on the Hirsch length $h(G)$ of $G$. When $h(G)=0$
then $G$ is finite and we can simply take $L=1$. Suppose that the Lemma has
been proved for all groups of Hirsch length less than $h>0$. Consider a virtually
polycyclic group $G$ with $h(G)=h$. Since $G$ is infinite it has an infinite normal abelian subgroup $N$. Let $\overline G= G/N^p$. Clearly $h (\overline
G) <h$ and so $\overline G$ has a $p$-good subgroup, say $\overline L = L/N^p$ for some normal subgroup $L$ containing $N^p$.

Then $L$ is $p$-good for $G$: Suppose $H \leq G$ with $HL=G$. Then $\overline
H \overline L= \overline G$ where $\overline H= HN^p/N^p$ and so the index $[G:HN^p]$
is finite and coprime to $p$. Therefore $HN=HN^p$ (since $N/N^p$ is a power of $p$ while $[HN:HN^p]$ must be coprime to $p$) and so $H_1 N^p =N$ where $H_1= H \cap N$. The last equality implies
that $[N:H_1]$ is finite and coprime to $p$ and hence so is $[G:H]=[G:NH][N:H_1]$.

Notice that the $p$-good subgroup $L$ we found contains $N^p$.
\end{proof}
Lemma \ref{pgood} raises the following natural
\begin{problem} Which groups $G$ possess a subgroup $N$ of finite index such
that $NH=G$ for a subgroup $H<G$ implies that $[G:H] < \infty$?
\end{problem}

The purpose of the following three Lemmas is to ensure the existence of suitable
word $w$, subgroup $V$ and elements $g_1, \ldots g_k$ in a virtually metabelian
and polycyclic group $G$ which meet the conditions of Theorem \ref{lifts}.
\begin{lemma}
\label{subgroups}
Let $G$ be a virtually metabelian and polycyclic group.
Then there exist normal subgroups $G_0 \rhd  V$ of $G$ such that
\begin{enumerate}
\item $G/G_0$ is finite, $V$ is a torsion free abelian group,
\item if $H$ is a subgroup of $G$ such that $HG_0=G$ then $H$ is of finite index in $G$,
\item $G_0/V$ is a nilpotent group which acts commutatively on $V$, i.e., $G_0/ C_{G_0}(V)$ is abelian, and

\item $\mathbb{Q} \otimes V$ is a perfect $\mathbb{Q}[G_0/V]$ module, i.e., $(G_0-1) \cdot V$ has finite
index in $V$.
\end{enumerate}
\end{lemma}
\begin{proof}
Let $A <B$ be normal subgroups of $G$ such that $A$ and $B/A$ are torsion
free abelian and $[G:B]$ is finite. Let $L$ be a $p$-good subgroup of $G$
for some prime $p$ and take $G_0=B \cap L$. Then $G_0$ also a $p$-good subgroup of $G$ and item 2 follows.

Let $W= A \cap G_0$. We have that $W$ and $G_0/W$ are torsion free abelian
groups and $W$ is a module for $\bar G_0 =G_0/W$.
Consider the chain of submodules $W\geq  (G_0-1)W \geq (G_0-1)^2 W \geq \cdots
$. This is a chain of subgroups of the finitely generated abelian group $W$
, so let $V$ be the first module in that series such that $(G_0-1)V$ has finite index
in $V$. Clearly $G_0/V$ is a nilpotent group (since $G_0$ acts
nilpotently on $W/V$). Item 4 is clear since $[V:(G_0-1)V]$ is finite while
item 3 follows since $C_{G_0}(V) \geq W$ and $G_0/W$ is abelian.
\end{proof}

\begin{lemma}
\label{action}
Let $\Gamma_0$ be a finitely generated torsion free abelian group and
let $V$ be a finitely generated torsion free $\mathbb Z [\Gamma_0]$ module such
that $V_\Q=V \otimes \mathbb Q $ is a perfect $\Q [\Gamma_0]$-module.
Then there exists an integer $N$ such that for any subgroup $\Gamma < \Gamma_0$
of index co-prime to $N$ we have that $V_\Q$ is a perfect $\Q[\Gamma]$-module.

Further when this happens then we can find an integer $M \in \mathbb N$ 
(depending on $\Gamma$),
integers $s_i$ and group elements $h_i \in  \Gamma$
such that
$$
\sum s_i (1-h_i) = M \cdot id
$$
as operators on $V$.
\end{lemma}
\begin{proof}
Let $\chi$ be a irreducible character (over $\C$) of $\Gamma_0$. We will call $\chi$ a
character of finite order if all values of $\chi$ are roots of $1$, in this case the
order of $\chi$ is the least integer $n$ such that all its values are $n$-th roots of $1$.

Since $V$ is a perfect $\Q[\Gamma_0]$ module it does not contain a trivial submodule.
Let $N$ be the gcd of all orders of irreducible characters which appears in $\C V$.
If $\Gamma\leq \Gamma_0$ is a subgroup of index co-prime to $N$ then the restriction of any
irreducible characters in $V$ to $\Gamma$ is non trivial. Therefore $\C V$ is a
perfect $\C[\Gamma]$ module, which implies that $V_\Q=\Q \otimes V$ is perfect $\Q[\Gamma]$ module.

For the second part, let $I$ be the augmentation ideal of $\Q \Gamma$. Since $V_\Q= IV_\Q$ and $V_\Q$
is a finite dimensional vector space over $\Q$ we have by Nakayama's lemma that
$id+T$ annihilates $V_\Q$ for some $T \in I$.

Expressing $T$ in the basis of $I$ and clearing the common denominator
$M$ of the rational coefficients gives the integers $s_i$ and
the elements $h_i \in G$.
\end{proof}

\begin{lemma}
\label{word}
Let $\Gamma$ be a group with an abelian normal subgroup $V$ such that $[\Gamma',V]=1$
and $\Gamma/V$ is nilpotent of class $d$.
Let $g_1,\dots, g_{k-1} \in \Gamma$.

For any integers $t \in \mathbb N, s_i$ and group elements $h_i \in \langle g_1,\dots, g_{k-1} \rangle$ \ ($i=1,\ldots , t$)
 there exists a word $w$ on $x_1,\dots,x_{k}$ such that
\begin{itemize}
\item $w(g_1,\dots,g_{k-1},g) \in V$ for any $g \in \Gamma$
\item The action of Fox derivative $\frac{\partial w}{\partial x_{k}}(g_1,
\ldots ,g_{k-1},g)$ on $ x \in V$ is given by
$$\displaystyle
x \mapsto \left( \sum_{i=1}^t s_i (1-h_i)\right)^d \cdot x.
$$
\end{itemize}
\end{lemma}
\begin{proof}
Consider the word $w'(g_1,\dots,g_{k-1},g)=\prod [g,h_i]^{s_i}$
where $h_i$ are expressed as words on $g_1,\dots,g_{k-1}$.
A direct computation gives that the Fox derivative 
$\frac{\partial w}{\partial x_k}$ at $(g_1, \ldots,
g_{k-1},g)$ with respect to the last variable acts on $V$ as multiplication by 
$\sum s_i (1-h_i) \in \Z [\Gamma/V]$
(use that $[ag, h_i]=  a [g,h_i] ({}^{h_i}a)^{-1}={}^{(1-h_i)}a \cdot [g,h_i]$
and each $[g,h_i]$ acts trivially on $V$).

Iterating the map $g \to w'(g_1,\dots,g_{k-1},g)$ \  $d$ times gives a
word $w$:
$$
w(x_1,\dots,x_k)=w'\bigg(x_1,\dots,x_{k-1}, w'\big(x_1,\dots,x_{k-1},w( \cdots w'(x_1,\dots,x_{k-1},x_{k})\dots\big)\bigg).
$$
The Fox derivative of $w$ with respect to $x_{k}$ is
$\big( \sum s_i (1-h_i)\big)^d$,
because substitution of words corresponds to multiplication of Fox derivatives.
The word $w$ always evaluates to one on $\Gamma/V$, because the group $\Gamma/V$ is nilpotent of class $d$.
\end{proof}

We now have all the ingredients to prove Theorem \ref{main}. It will follow
from the corresponding result for metabelian groups:

\begin{theorem}
\label{met}
Let $G$ be a virtually metabelian polycyclic group.
Then
$$
d(\widehat{G}) \leq d(G) \leq d(\widehat{G})+1.
$$
\end{theorem}
\begin{proof}
Let $k=d(\widehat{G})+1$. Let $G_0$ and $V$ are the subgroups provided by
 Lemma~\ref{subgroups}. Now Lemma~\ref{action} applied to the group $G_0/[G_0,G_0]V$
acting on $V$,
gives us an integer $N$ such than any subgroup of index co-prime to $N$ in $G_0/V$
acts perfectly on $V$ (as a rational module).

By Lemma \ref{vabelian} there exists a generating set
$S=\{\gamma_1,\ldots,\gamma_k\}$ of $G/V$ such that
$\gamma_k\in G_0$ and the subgroup $\Gamma=\langle \gamma_1 , \ldots ,\gamma_{k-1}
\rangle  V \cap G_0$ has index $[G_0:\Gamma]$ co-prime to $N$. Therefore $V$ is a perfect
rational $\Q [\Gamma/V]$-module and for some integer $M$ and element $T$
in the augmentation ideal of $\mathbb Z [\Gamma/V]$ we have that $T$ acts
on $V$ as multiplication by $M$.

For each prime $p$ pick a $p$-good subgroup $G_p$ of $G$ containing $V^p$. (In fact by replacing $G_p$ with a normal subgroup of $G$ of finite index
we may even assume $G_p \cap V=V^p$.)

Now apply Lemma~\ref{word} to $\Gamma$ and $V$ with $\gamma_i=g_i$, ($i=1,\ldots ,k-1$),
$g=\gamma_k$ and $s_i \in \mathbb N$, $h_i
\in \langle \gamma_1, \ldots ,\gamma_{k-1}\rangle$ chosen so that $T \equiv \sum_i s_i(1-h_i)$ in $\mathbb Z[\Gamma/V]$. We conclude that there is a word $w(x_1,
\ldots , x_k)$ such that $w(\underline \gamma) \in V$ and 
$\frac{\partial w}{\partial x_k}(\underline \gamma)$ acts on $V$ as multiplication
by $M^d$ where $d$ is the nilpotency class of $\Gamma/V$.

We can now apply Theorem~\ref{lifts} to $G$ with these choices of $w,V, G_p$ and $\gamma_j$. The conditions 1,2 and 3 are satisfied by the construction of
the subgroups $G_0$ and $V$, the word $w$ and the definition of the number $k$. So by Theorem \ref{lifts} we can find lifts $a_i \in \gamma_i V$ such that
$G=\langle a_1, \ldots ,a_k \rangle$. Theorem \ref{met} is proved.
\end{proof}

\begin{remark} 
\label{equality}
If we have that $d(G/V) < d(\widehat G)$ then the argument above gives that
$d(G) = d(\widehat G)$.
\end{remark}
\bigskip

\begin{proof}[Proof of Theorem \ref{main}] In general a virtually polycyclic group $G$ is virtually nilpotent by abelian,
i.e., it has normal subgroups $G_1>G_2$ such that $G/G_1$ is finite, $G_1/G_2$
is abelian while $G_2$ is nilpotent. (See Theorem 2, Chapter 2 in \cite{Danbook}).

Now every group which generates $H=G/G_2'$
generates $G$ and so we have $d(G)=d(H),\ d(\widehat G)=d(\widehat H)$. Thus
Theorem \ref{main} becomes a corollary of Theorem \ref{met}. Moreover its
proof gives an efficient algorithm for generating a polycyclic group $G$
with $d(\widehat G)+1$ elements, even with $d(\widehat G)$ elements if the condition
of Remark \ref{equality} holds. \end{proof}

\subsection*{Acknowledgement}
We thank B.~Petrenko for his preprint~\cite{Petrenko},
which contains an argument on the generation of rings which led us to
Theorem~\ref{lifts}. Similar result for modules can be found in~\cite{Swan}.
The first author thanks Imperial College, London
for the hospitality during his visit.

\end{document}